\documentclass[12pt,leqno]{article}
\usepackage{amssymb,amsmath,amsthm,amscd,latexsym,verbatim,enumerate}
\usepackage[pagebackref]{hyperref}
\usepackage{amsrefs}
\usepackage{float}
\usepackage[mathscr]{eucal}
\usepackage{color}
\hypersetup{citecolor=red, linkcolor=blue, colorlinks=true}

\newtheorem{thm}{Theorem}[section]
\newtheorem{prop}[thm]{Proposition}

\newtheorem{cor}[thm]{Corollary}

\newtheorem{defi}{Definition}
\newtheorem{remk}{Remark}
\newcommand{\pf}{{\bf Proof. \ }}

\renewcommand{\geq}{\geqslant}
\renewcommand{\leq}{\leqslant}
\renewcommand{\ge}{\geqslant}
\renewcommand{\le}{\leqslant}

\newcommand{\F}{\mathbb{F}}

\newcommand{\Ker}{\textup{Ker}}
\newcommand{\Synd}{\textup{Synd}}

\def\1v{\mathbf{1}}
\def\0v{\mathbf{0}}

\def\GL{{\rm GL}}

\begin{document}
\title{ Twisted Centralizer Codes}
\author{Adel Alahmadi\footnote{Mathematics Department, King Abdulaziz University, Jeddah, Saudi Arabia.}\\
S.\,P. Glasby\footnote{Center for Mathematics of Symmetry and Computation, University of Western Australia, 35 Stirling Highway, WA 6009, Australia. Email: \texttt{Stephen.Glasby@uwa.edu.au, Cheryl.Praeger@uwa.edu.au}}\;\footnote{The Department of Mathematics, University of Canberra, ACT 2601, Australia.}\\
Cheryl E. Praeger$^\dagger{}^*$\\
Patrick Sol\'e\footnote{
CNRS/LAGA, University of Paris 8, 2 rue de la Libert\'e, 93 526 Saint-Denis, France.
Email: \texttt{sole@enst.fr}}\\
Bahattin Yildiz\footnote{University of Chester, UK. Email: \texttt{bahattinyildiz@gmail.com}}}
\date{\today}
\maketitle

\let\thefootnote\relax

\begin{abstract}
Given an $n\times n$ matrix $A$ over a field $F$ and a scalar $a\in F$,
we consider the linear codes
$C(A,a):=\{B\in F^{n\times n}\mid \,AB=aBA\}$ of length~$n^2$. We call $C(A,a)$ a
\emph{twisted centralizer code}.
We investigate properties of these codes including their dimensions, minimum distances, parity-check matrices, syndromes, and automorphism groups. 
The minimal distance of a centralizer code (when $a=1$) is at most~$n$,
however for $a\ne0,1$ the minimal distance can be much larger, as large as~$n^2$.
\end{abstract}

\section{Introduction}
Denote the $n\times n$ matrices over a field $F$ by $F^{n\times n}$.
Fix a matrix $A\in F^{n\times n}$ and a scalar $a\in F$. As we are motivated by applications to coding theory we focus on
the case where $F$ is a finite field $\F_q$ of order $q.$ The \emph{centralizer of $A$, twisted by $a,$} is defined to be
\begin{equation}
C(A,a):=\{B\in F^{n\times n}\mid \,AB=aBA\}.
\end{equation}
Clearly $C(A,a)$ is an $F$-linear subspace of the vector space $F^{n\times n}$. 
Note that $C(A,0)$ is the right-annihilator of $A$. We shall use the notation $C(A)$ instead of $C(A,1)$ when $a=1$, and note that $C(A)$ is simply the centralizer of $A$.
The subspace $C(A,a)$ of $F^{n\times n}$ is viewed as code: we view a codeword $B\in F^{n\times n}$ as column vector $[B]$ of length $n^2$, by reading the matrix $B$ column-by-column.
The case $a=1$ was considered in \cite{centralizer}. In the present paper, we extend and sometimes correct the results of \cite{centralizer}.
In particular the incorrect \cite[Theorem~ 2.4]{centralizer} is corrected and generalized for this larger class of codes in \cite{GP} (see Theorem \ref{low}), 
and we exploit this result in several ways in Section~\ref{sect2.2}. 

\begin{defi}{\rm 
For any $n\times n$ matrix $A\in F^{n\times n}$ and any scalar $a\in F,$ the subspace $C(A,a)$ formed above is called the {\it centralizer code} obtained from $A$
and \emph{twisted by}~$a.$
}
\end{defi}
In a sense $A$ serves as a parity-check matrix, because $B$ lies in $C(A,a)$ precisely when $AB-aBA = 0$.
More concretely, in the following result we show that a certain $n^2\times n^2$ matrix $H$ related to $A$ is a parity check matrix 
in the sense that $B\in C(A,a)$ if and only if $H[B]=0$, where $[B]$ is the $n^2$-dimensional column vector above corresponding to $B$.

{\prop \label{prop1.1} A parity-check matrix for $C(A,a)$ is given by 
$$
H=I_n\otimes A-a(A^t\otimes I_n),
$$ 
where $\otimes$ denotes the Kronecker product,
 and $A^t$ the transpose of the matrix $A.$}

\pf
This follows from the proof of \cite[Theorem 27.5.1, p.124]{P} with $A=A$, $B=a A$, $C=0.$ Also, a direct proof (for row vectors) 
is given in \cite[Lemma 3.2]{GP}.
\qed

\medskip
The following simple observations involving $C(A,a)$ will be used later.

{\thm \label{nite}
 Suppose $a,a'\in F$ and $A \in F^{n \times n}$ is a matrix. Then the following are true:
\par {\bf a)} $A \in C(A,a)$ if and only if $a = 1$ or $A^2 = 0$.
\par {\bf b)} If $B\in C(A, a)$ and $B' \in C(A,a')$ then $BB'$ and $B'B$ both lie in $C(A,aa')$.
\par {\bf c)} For $a \neq 0$, $B \in C(A, a) \Leftrightarrow A \in C(B,a^{-1})$.
\par {\bf d)} For $A\neq 0_{n\times n}$, we have $I_n \in C(A,a) \Leftrightarrow a=1$.
\par {\bf e)} For $a\ne 0$,  $B\in C(A,a)  \Leftrightarrow B^t\in C(A^t,a^{-1})$.
}

\pf We only prove (b) since the other parts are simple observations.
Starting from $AB=a BA$ and $AB'=a'B'A,$ successive substitutions give
$$A(BB')=(AB)B'=aBAB'=aa'(BB')A. $$
\qed

\newpage

The problems about $C(A,a)$ that arise naturally for given $A$ and $a$ include:
\begin{itemize}\itemsep0em
\item computing its dimension ($k$);
\item deriving decoding/encoding algorithms and bounding the minimum distance ($d$);
\item determining its automorphism group.
\end{itemize}

The paper is organized as follows. Sections 2, 3, 4 successively tackle the above three problems.
Section 5 is dedicated to concrete examples of codes and Section 6 contains some concluding remarks and open problems.  We say that $C(A,a)$ has \emph{parameters}
$[n^2,k,d]$ where $A\in F^{n\times n}$, $k=\dim(C(A,a))$, and $d$ is the minimal
(Hamming) weight of a nonzero vector. Note that $I_n\in C(A,1)$, so the minimal distance of a centralizer code is at most~$n$.

\section{Dimension}
\subsection{Basic Bounds}\label{sect2.1}

{\prop \label{L:ExistInvertible}
If $A\in F^{n\times n}$, $a\in F$, and $C(A,a)$ contains an invertible matrix,~then
\[
  \dim C(A,a)=\dim C(A).
\]}

\pf
The result is true for $A=0$ as $C(0,a)=F^{n\times n}=C(0)$.
Suppose now that $A\ne0$ and that $B\in C(A,a)$ is invertible.
Then $B\in C(A,a)$ implies that $A=aBAB^{-1}$, and since $A\ne0$ we must have $a\ne0$.  
The linear map
$\phi_B\colon C(A)\to C(A,a)$ with $X^{\phi_B}= XB$
is injective so
$\dim C(A)\le\dim C(A,a)$. 
Then
$B^{-1}\in C(A,a^{-1})$ and $YB^{-1}\in C(A)$ for all
$Y\in C(A,a)$. 
Hence the map $\psi_B\colon C(A,a)\to C(A)$ with
$Y^{\psi_B}=YB^{-1}$ is the inverse of $\phi_B$ above.
Therefore $\psi_B$ is an isomorphism and
$\dim C(A,a)=\dim C(A)$ as claimed.
\qed

\medskip


The next result shows that for every $n$, with mild assumptions on $A$ and $a$, the dimension of $C(A,a)$ is bounded above.

{\prop If $0\neq A \in F^{n\times n}$ and $a\ne 1,$ then $\dim(C(A,a))\le n^2-1.$}

\pf 
If $\dim(C(A,a))=n^2,$ then every matrix $B$ satisfies the relation $AB=aBA.$ 
However if $B=I$ then $A=aA$ which does not hold.
\qed

In fact this bound can be improved to $\dim(C(A,a))\le n^2-n$, see \cite[Corollary 6.7]{GP}.

\subsection{Spectral bounds on the dimension of \texorpdfstring{$C(A,a)$}{} }\label{sect2.2}
 In this section we treat matrices over the finite field $\F_q$.  Some spectral notation is in order.
For $B\in \F_q^{n \times n},$ let $F$ denote a splitting field for the characteristic polynomial of $B$, 
and denote by $S(B)\subseteq F$ the set of its eigenvalues in $F$.
 Let $K(B,\lambda)$ denote the dimension over $F$ of the eigenspace (in $F^n$) corresponding to $\lambda \in S(B).$
Denote by $M(B,\lambda)$ the multiplicity of $\lambda$ as a root of the characteristic polynomial of $B.$
From linear algebra we know that $K(B,\lambda)\le M(B,\lambda)$.

{\thm {\rm \cite[Theorem 4.7]{GP}} \label{low} Let $a\in \F_q$ and $A\in \F_q^{n \times n}$. Then 

$$
\sum_{ \mu\in S(A)} K(A,a\mu)K(A^t,\mu)\leq \dim_{\F_{q}}(C(A,a))\leq \sum_{ \mu\in S(A)} M(A,a\mu)M(A^t,\mu) .
$$
}

An important consequence of the lower bound is the following, where $\Ker(A)$ denotes the null space of $A$.

{\cor $$\dim C(A,a) \ge (\dim \Ker(A) )^2.$$}

\pf
In the formula in Theorem~\ref{low}, we may bound the sum below by the term indexed by $\mu=0.$ Note that $\dim \Ker(A)=\dim \Ker(A^t).$
Alternatively we may invoke Theorem \ref{product} below.
\qed

\medskip
The next Corollary explains why many matrices  $A$ yield codes $C(A,a)=\{0\}.$

{\cor\label{C:0} If there is no $\lambda \in S(A)$ such that $\lambda a \in S(A)$ then $\dim(C(A,a))=0.$ }

\pf
In the upper bound in Theorem~\ref{low}, all the summands are zero.
\qed

The bounds are most useful when $A$ is a combinatorial matrix with a known spectrum. Recall that an Hadamard matrix of order $n$ is a $\{\pm 1\}$ valued matrix $H$ satisfying
$HH^t=nI$, see~\cites{MS,W}. 
The Kronecker product $H_{2^k}$ of $k$ copies of
$\left[\begin{smallmatrix}1&1\\1&-1\end{smallmatrix}\right]$ are
examples of degree $2^k$.

\begin{prop} \label{prop2.7}
Suppose $A\in\F_q^{n\times n}$ where $q$  is a power of an odd prime $p$, and
$p\nmid n$. If $A$ is a symmetric Hadamard matrix with trace zero,
  then  $\dim(C(A,a))$ equals $\frac{n^2}{2}$ if $a=\pm1$, and $0$ otherwise.
\end{prop}

\pf 
Note that $p=\textup{char}(\F_q)\nmid n$.  
By definition, $AA^t=nI,$ and by hypothesis $A=A^t,$ so that $A^2=nI.$
Thus the eigenvalues of $A$ are $\pm\lambda$ where $\lambda^2=n\ne 0$ in $\F_q$.
Hence $\lambda$ lies in $\F_{q^2}^\times$, and the Jordan form of $A$ must be
$\lambda I_k\oplus (-\lambda)I_{n-k}$
for some $k$ with $0<k<n$. However, $0=\textup{Trace}(A)=k\lambda+(n-k)(-\lambda)$ implies $k=n/2$. It follows that $K(A,\pm\lambda)=M(A,\pm\lambda)=n/2$.
Thus if $a=\pm1$ the upper and lower bounds of Theorem~\ref{low} coincide,
and $\dim(C(A,a))=2(n/2)^2=n^2/2$. If $a\ne\pm1$, then $\dim(C(A,a))=0$
by Corollary~\ref{C:0}.
%
\qed

\medskip\noindent
{\bf Example 1:} Let $F$ be one of the fields $\F_3$ or $\F_5$. Then, taking a Sylvester-type Hadamard matrix 
$H_4$ of order $n=4$ yields an isodual$^*$\footnote{$^*$An \emph{isodual} code is
one that is monomially equivalent to its dual.} code $C(H_4,-1)$ over $F$ with parameters $[16,8,4].$
The next order, that is, $H_4\oplus H_4$ gives a $[64,32,8]$ code. Note that $8$ is not a 
quadratic residue modulo $3$ or modulo $5$.


\section{Encoding-Decoding}

Our approach to encoding-decoding procedures is similar to the case of ordinary centralizer 
codes $C(A)$ discussed in \cite{centralizer}. 
The codes $C(A,a)$ retain the advantage of efficient
syndrome computation of the ordinary centralizer codes. An important difference is the much higher error 
correction capability of twisted codes with respect to the highly restricted capacity of the centralizer 
codes $C(A)$. If $C(A,a)$ has dimension $k$, then the \emph{information rate} is $k/n^2$, and we can give a 
procedure for encoding and decoding. As an $\F_q$-vector space, $C(A,a)$ has a basis consisting of $k$ matrices, 
which we denote by $\{A_1,A_2, \dots, A_k\}$.
We encode a given information vector (or message) $(a_1,a_2,\dots, a_k) \in \F_q^{k}$ as a codeword in $C(A,a)$, namely
$$
a_1A_1+a_2A_2+ \dots + a_kA_k.
$$
The decoding can be done by reversing the above procedure. So, to find the message that a matrix $B\in C(A,a)$ represents, 
all we need is to find the coordinate vector for $B$ in the basis $\{A_1, A_2, \dots, A_k\}$ of $C(A,a)$.

Note that $C(A,a)$ is an additive subgroup of $\F_q^{n\times n}$, and hence partitions $\F_q^{n\times n}$ into its additive cosets. 
Since a matrix  $B\in \F_q^{n\times n}$ lies in $C(A,a)$ if and only if $AB-aBA=0$, we can use $A$ itself
(instead of the $n^2 \times n^2$ matrix $H$ in Proposition~\ref{prop1.1}) as a type of 
parity-check matrix to obtain `syndromes' which are determined by the cosets of $C(A,a)$.
 Thus we make the following definition:

\begin{defi}{\rm 
Let $B \in \F_q^{n\times n}$ be any matrix over $\F_q$.
The {\it syndrome} of $B$ in $C(A,a)$ is defined as
$$
\Synd_{A,a}(B) = AB-aBA.
$$
}
\end{defi}
Thus $\Synd_{A,a}$ is an $F$-linear map $F^{n\times n}\to F^{n\times n}$ and
$B\in C(A,a) \Leftrightarrow \Synd_{A,a}(B) = 0.$
The following theorem suggests that this definition of the syndrome might help us in an error-correction scheme.

{\thm \label{thm:synd} Let $B_1, B_2 \in \F_q^{n\times n}$. Then $\Synd_{A,a}(B_1) = \Synd_{A,a}(B_2)$ if and only if $B_1$ and $B_2$ are 
in the same additive coset of $C(A,a)$.  }

\pf
The proof is similar to the proof of  \cite[Theorem 3.1]{centralizer} and is therefore omitted.
\qed

Testing whether a matrix $B$ lies in $C(A,a)$ is the same as checking
whether $\Synd(B)=0$. (We henceforth drop the subscripts on $\Synd$.)
Multiplying two $n\times n$ matrices has computational complexity of
$O(n^m)$ field operations where $m\le 2.373$, see the survey in~\cite{V}.
Thus testing whether $B$ lies in $C(A,a)$ using syndromes has the
same complexity. Alternatively, one could use the $n^2\times n^2$ parity
check matrix $H$ of Proposition~\ref{prop1.1}.
Since multiplying a vector in $\F_q^{n^2}$ by $H$
(via the na\"{\i}ve algorithm)
requires $O(n^4)$ field operations, the syndrome method is
computationally advantageous.



We show in Subsection~\ref{sub:mind} that twisted centralizer codes
have higher minimum distances, and hence higher error correction capability,
than centralizer codes.

\subsection{Bounds on the minimum distance}\label{sub:mind}
The distance $d$ of a nonzero linear code is the minimal Hamming weight (number of nonzero coordinates) of a nonzero vector in the code. Thus the distance
of $C(A,a)$ is at most $n^2$. 

Let $J_n$ denote the $n\times n$ matrix with all entries equal to~1. 
The following theorem shows that there exist codes $C(A,a)$ whose
minimal distance $d$ is $n^2$, which is as large as possible.

{\thm \label{thm:dd}Suppose $J_n+I_n\in F^{n\times n}$ where
$\textup{char}(F)$ divides $n+1$. If $a\ne 0,1$, then the twisted centralizer
code $C(J_n+I_n,a)$ equals $\langle J_n\rangle$ and has parameters $[n^2,1,n^2]$.}

\pf
Let $A=J_n+I_n$, and let $p=\textup{char}(F)>0$.
It follows from $J_n^2=nJ_n$ 
that $J_n^p=n^{p-1}J_n=(-1)^{p-1}J_n=J_n$. Hence $A^p=(J_n+I_n)^p=J_n^p+I_n^p=A$.
Let $u\in F^n$ be the column vector with all entries 1. Since
$J_nu=nu$ and $Au=(n+1)u=0$, we see that $\det(A)=0$.
The null space of $J_n$ has dimension~$n-1$, and so the same is true
of the 1-eigenspace of $A$. Thus the Jordan form of $A$ is
the diagonal matrix $D=\textup{diag}(0,1,\dots,1)$.  A direct calculation
using the fact that $a\ne0,1$ shows that $C(D,a)=\langle E_{11}\rangle$
is 1-dimensional.
Since~$A$ and $D$ are conjugate, the same is true for $\dim C(A,a)$.
However, $J_nA=AJ_n=(n+1)J_n=0$ and hence $C(A,a)=\langle J_n\rangle$.
Therefore each nonzero element of $C(A,a)$ is $bJ_n$ for some nonzero $b$,
and hence $C(A,a)$ has minimum distance $n^2$.
In summary, $C(A,a)$ has parameters $[n^2,1,n^2]$ as claimed.
\qed

A code with parameters $[N,k,d]$ must have $k+d\le N+1$. Hence if $N=d=n^2$,
we conclude that $k=1$ as happens in Theorem~\ref{thm:dd}.

Experimental evidence using the programs~\cite{Glasby} suggests that
for any field other than $\F_2$, and for any scalar other than 0 or 1,
there exists some matrix
$A\in F^{n\times n}$ such that $C(A,a)$ has parameters $[n^2,1,n^2]$. Proving this claim for all $n\ge2$ appears to be difficult
without first guessing the \emph{form} of suitable matrices~$A$, and even then
the computations can depend in a complicated way on the field $F$. Let
$A_n=E_{n,1}-E_{n,n}+\sum_{i=1}^{n-1}E_{i,i}-E_{i,i+1}$, and let
$B_n=J_n-2\sum_{i=1}^{n}E_{i,n}$. When $n=4$ these matrices are
\[
  A_4=\begin{pmatrix}1&-1&0&0\\0&1&-1&0\\0&0&1&-1\\1&0&0&-1\end{pmatrix},
  \qquad\textup{and}\qquad
  B_4=\begin{pmatrix}1&1&1&-1\\1&1&1&-1\\1&1&1&-1\\1&1&1&-1\end{pmatrix}.
\]
If the characteristic of $F$ is zero, then a variant of the proof of
Theorem~\ref{thm:dd} 
shows that $C(A_n,-1)=\langle B_n\rangle$. However, our focus is on
\emph{finite} fields, and here one can only prove this when
$\textup{char}(F)$ is `sufficiently large'. The idea is to solve the system
$A_nX+XA_n=0$ for $X\in\mathbb{Z}^{n\times n}$, thereby determining `bad' primes
which will increase the rank of the solution space.
For example, $C(A_n,-1)=\langle B_n\rangle$ holds when $n=3$ provided
$\textup{char}(F)\ne2$, and when $n=30$ provided $\textup{char}(F)$ is not
one of 27 `bad' characteristics. The bad characteristics have
$k=\dim C(A_n,-1)>1$, and hence $d<n^2$. They have much larger information
rates $k/n^2$, but may correct fewer errors.

\subsection{An upper bound}
An $\F_q$-linear code $C$ is said to have \emph{parameters} $[N, K, d]$
if $C$ is a $K$-dimensional subspace of $\F_q^N$ (consisting of column vectors),  
$d$ is the minimum distance of $C$, and $N$ is called the \emph{length} of $C$. 
Sometimes we omit $d$ and say that $C$ is an $[N,K]$ code over $\F_q$.
If $D,E$ are two $F$-linear codes of the same length ~$N$, we write their {\em product code} as

$$ D\otimes E=\{ uv^t|\, u \in D,\, v \in E\}\subseteq F^{N\times N}.$$
{\thm \label{product} For all $a\in \F_q,$ the code $C(A,a)$ contains the product code $\Ker(A)\otimes \Ker(A^t).$ If $\Ker(A),\Ker(A^t)$ have respective parameters $[n,k,d]$ and
$[n,k',d'],$ then $C(A,a)$ has parameters $[n^2,K,D]$ with $K\ge kk'$ and $D\le dd'.$}

\pf
If $u\in \Ker(A),$ and $v\in \Ker(A^t),$ (both column vectors) then $B=uv^t \in C(A,a),$ since $AB=(Au)v^t=0,$ and $BA=u(A^tv)^t=0.$ The second statement follows by
standard properties of product codes \cite{MS}.
\qed


\medskip
This result leads to a general upper bound on the minimum distance $d(A,a)$ of $C(A,a).$  
Denote by $\Delta_q(N,K)$ the largest minimum distance of all $[N,K]$ codes over $\F_q.$

{\cor For all $a\in\F_q,$ and $A\in\F_q^{n\times n}$, 
$$
d(A,a) \le (\Delta_q(n,k_0))^2, \quad\mbox{where}\quad k_0=\dim_q(\Ker(A)).
$$ }

\pf Both $\Ker(A)$ and $\Ker(A^t)$ have the same dimension $k_0.$ The result follows by the second statement of Theorem \ref{product}. \qed

\medskip
Theorem~\ref{product} also gives the following lower bound for the minimum distance $d(A,a)$ of $C(A,a)$ for rank 1 matrices $A$.
 
{\cor \label{cor-rank1} If $A$ has rank one then  $d(A,a)$ is at most $4.$}

\pf
Since $A$ has rank 1, both $\Ker(A)$ and $\Ker(A^t)$ have dimension $n-1$. 
Hence either $\Ker(A)$ contains a weight~1 vector, or two distinct weight 1 vectors not in $\Ker(A)$ will differ by a 
weight 2 vector, which must lie in $\Ker(A)$. Thus $\Ker(A)$ contains a nonzero vector of weight at most $2$, and hence
(when regarded as a code in $\F_q^n$) it has minimum distance at most $2$. The same is true for $\Ker(A^t)$, 
and the result now follows from the second statement of Theorem \ref{product}.
\qed

\subsection{A lower bound: asymptotics}

Recall the $q-$ary entropy function defined for $0<x< \frac{q-1}{q}$ by 
$$ 
H_q(x)=x\log_q(q-1)-x\log_q(x)-(1-x)\log_q(1-x).
$$
This quantity is instrumental in the estimation of the volume of high-dimensional Hamming balls when the base field is $\F_q.$
The result we are using is that the volume of the Hamming ball of radius $xn$ in $\F_q^n$ is, up to subexponential terms,  
$q^{nH_q(x)},$ when $0<x<1$ and $n$ goes to infinity
\cite[Lemma 2.10.3]{HP}.

{\thm For $a \neq 0,1$, and for $n \rightarrow \infty,$  
there are codes $C(A,a)$ with minimum distance at least $ \frac{n}{\log n}.$ }

\pf
Let $B\in\F_q^{n\times n}$ be nonzero and have weight less than $\frac{n}{\log n}$. Then $B\in C(A,a)$ if and only if $A\in C(B,a^{-1})$,
by Theorem~\ref{nite}(c). Thus the number of codes $C(A,a)$ containing $B$ is $|C(B,a^{-1})|$, and by   \cite[Prop. 6.6]{GP},
this is at most $X := q^{n^2-2n+2}$.

As we mentioned above, by \cite[Lemma 2.10.3]{HP}, the number of such matrices $B$ is $q^{n^2H_q(Y)}$ where $Yn^2= \frac{n}{\log n}$, that is, 
$Y=\frac{1}{n\log n}$. Thus the total number of codes $C(A,a)$ which contain at least one nonzero matrix $B$ of weight less than $\frac{n}{\log n}$ is at most 
$X q^{n^2H_q(Y)}=q^{n^2-2n+2+ n^2H_q(Y)}$. The result will follow if we can show that this quantity is less than the total number of these codes $C(A,a)$
which are nonzero. By \cite[Theorem 5.2]{GP}, the total number of these nonzero codes is at least $\frac{1}{q}\times q^{n^2} = q^{n^2-1}$.

Now $H_q(x)=-x\log_q x +O(x),$ for small $x$ (see, for example, \cite[Proposition 3.3.6]{GRS}), and hence
for $n\rightarrow \infty,$ we have
$$
n^2 H_q(Y) = n^2H_q(\frac{1}{n\log n})\sim -\frac{n}{\log n}\log_q \frac{1 }{n\log n}\sim  \frac{n}{\log q} . 
$$
Thus $n \sim n^2H_q(Y)\log q = \log (q^{n^2 H_q(Y)})$, and hence $q^{n^2 H_q(Y)} \sim q^n$.

Therefore, for large $n$, the number of codes $C(A,a)$ which contain at least one nonzero matrix $B$ of weight less than $\frac{n}{\log n}$ is at most 
$q^{n^2-2n+2+ n^2H_q(Y)} \sim q^{n^2-n+2}$, which is less than $q^{n^2-1}$.
 We can now conclude that there exist 
codes in the family that have minimum distance at least $\frac{n}{\log n}.$
\qed

%
%
%
%

\subsection{An example of error-correction for rank 1 matrices}\label{sub:rank1}
We discuss the twisted centralizer codes $C(A,a)$ for the case where $A$ is a rank 1 matrix in $\F_q^{n\times n}$. 
The dimensions of such codes were determined in \cite[Remark 2.10]{GP}, and all are of the form 
$(n-1)^2 + \delta$. Indeed if $a\ne 0, 1$, then $\delta = 1$ when $A^n=0$
and $\delta= 0$ otherwise.
However the minimum distance $d(A,a)$, and hence the error-correcting properties, are not so uniformly described.

For example, if $A=E_{11}$, the matrix with entry $1$ in the $(1,1)$-position and all other entries zero, then for any value of $a$, 
the code $C(E_{11},a)$ is easy to compute, and in particular for each of the $(n-1)^2$ pairs $(i,j)$ with $i>1$ and $j>1$, $C(E_{11},a)$ contains the weight~$1$ matrix $E_{ij}$ 
(with a single nonzero entry, namely an entry 1 in the $(i,j)$-position).  Thus $d(E_{11},a)=1$, which is unfortunately smaller than the upper bound given in Corollary~\ref{cor-rank1}.

We show in Theorem~\ref{thm:rank1} that the upper bound of Corollary~\ref{cor-rank1} is often achieved for a different family of rank 1 matrices, 
namely the matrices $J_n$ where $J_n$ has degree~$n$ and has all entries equal to $1$. This illustrates, 
in particular, that conjugating $C(A,a)$ by an element of $\GL(n,q)$ can change
the minimal distance of the code.

The `single errors' that may occur are the weight $1$ matrices $b E_{ij}$, where $b\in\F_q\setminus\{0\}$. We say that a code $C(A,a)$ \emph{corrects single errors}
if distinct single errors $B, B'$ give distinct syndromes $\Synd(B)$ and $\Synd(B')$. Our next result shows that the codes $C(J_n,a)$ can be used for single error correction.

{\thm \label{thm:rank1} Let $q\geq3$, $n\geq2$, and $a\in\F_q\setminus\{0,1\}$. Then the code $C(J_n,a)$ corrects single errors.
Moreover,
\begin{enumerate}
 \item[{\rm (a)}] if either $n\geq3$ or $q$ is odd, then $d(J_n,a)=4$ and, for example, $E_{11} - E_{12} - E_{21} + E_{22}$ is a minimum weight nonzero codeword; while 
\item[{\rm (b)}] if $n=2$ and $q$ is even (so $q\geq4$), then $d(J_n,a)=3$ and $aE_{11} + (a-1)E_{12} + E_{22}$ is a minimum weight nonzero codeword. 
\end{enumerate}
}

\pf
As discussed above, a single error is a matrix $b E_{ij}$, with $b\ne0$ and $1\leq i, j\leq n$. 
A simple computation shows that the syndrome $\Synd(bE_{ij}) = J_n(bE_{ij})-a(bE_{ij})J_n = b S(i,j)$, where $S(i,j)=\Synd(E_{ij})$ is the matrix with
$(k,\ell)$-entry described below
\[
 S(i,j)_{k\ell} = \begin{cases}
  0 & \mbox{if $k\ne i$ and $\ell\ne j$}\\                   
  1 & \mbox{if $k\ne i$ and $\ell= j$}\\                   
  -a & \mbox{if $k= i$ and $\ell\ne j$}\\                   
  1-a & \mbox{if $k= i$ and $\ell= j$.}
                    \end{cases} 
\]
In particular, since all of these syndromes are nonzero, it follows that $C(J_n,a)$ contains no weight one matrices, so $d(C(J_n,a))\geq2$.

Next we consider a weight two matrix $B=b_1E_{i_1j_1} + b_2E_{i_2j_2}$ with $(i_1,j_1)\ne (i_2,j_2)$, $i_1\le i_2$, $b_1\ne0$ and $b_2\ne0$.
Then $S:= \Synd(B) = b_1 S(i_1,j_1) + b_2 S(i_2,j_2)$. Suppose first that $i_1\ne i_2$ and $j_1\ne j_2$. Up to row and column permutations, the $\{i_1,i_2\}\times\{j_1,j_2\}$-submatrix of $S$ is 
\begin{equation}\label{eq:synd}
 \begin{pmatrix} b_1(1-a)&b_2-ab_1 \\  b_1-ab_2	&b_2(1-a) \\ \end{pmatrix}. 
\end{equation}
Since $a\ne 1$ we see $b_1(1-a)\ne0$, so $S\ne0$.
Next suppose that $i_1 = i_2$, so that we must have $j_1\ne j_2$ since $B$ has weight two.
For $i\ne i_1$, the entry $S_{i,j_1}=b_1\ne0$, so $S\ne0$. A similar argument shows that $S\ne 0$ if $j_1=j_2$. Thus $S\ne0$ for all weight two matrices $B$,
which implies that   $d(C(J_n,a))\geq3$.  
 
A straightforward computation shows that, if $B=E_{11} + E_{22} - E_{12} - E_{21}$, then $\Synd(B)=0$, and hence $B\in C(J_n,a)$. 
Thus it remains to consider whether $C(J_n,a)$ contains a weight three codeword $B$.   
Suppose then that $B=b_1E_{i_1j_1} + b_2E_{i_2j_2}  + b_3E_{i_3j_3}$ with distinct pairs $(i_k,j_k)$ and nonzero $b_k$, for $k=1,2,3$,
and let $S := \Synd(B) = \sum_{k=1}^3 b_k S(i_k,j_k)$. Suppose that $B\in C(J_n,a)$, or equivalently, that $S=0$.

First assume that  $i_1\ne i_2$ and $j_1\ne j_2$, and let $S'=\Synd(b_1E_{i_1j_1} + b_2E_{i_2j_2})$. 
The $\{i_1,i_2\}\times\{j_1,j_2\}$-submatrix of $S'$ is as in \eqref{eq:synd}.
As each of the diagonal entries of this restriction is nonzero, and as $S=S'+b_3S(i_3,j_3)=0$,
it follows that (i) either $i_3=i_1$ or $j_3=j_1$; and also (ii)  either $i_3=i_2$ or $j_3=j_2$. Considering transposes if 
necessary, we may therefore assume that $i_3=i_1$ and $j_3=j_2$. If $n\geq3$ then, for $i\not\in\{i_1,i_2\}$, the entry $S_{ij_1}=b_1\ne 0$, which is a contradiction. Thus $n=2$, and
\begin{equation*}
S =  \begin{pmatrix}
  b_1(1-a)-b_3a	&	b_2-ab_1 +b_3(1-a)\\
  b_1-ab_2	&	b_2(1-a) +b_3 \\
 \end{pmatrix} = 0. 
\end{equation*}
It follows that $b_1=b_2a$ (from the $(i_2,j_1)$-entry), $b_3 = - b_2 (1-a)$   (from the $(i_2,j_2)$-entry), and then from these equalities and the
$(i_1,j_1)$-entry we find 
\[
0 = b_1(1-a) - b_3 a = 2b_2a(1-a)   
\]
so that, since $b_2a(1-a)\ne0$, we must have $q$ even. The $(i_1,j_2)$-entry is $2b_2a(1-a)$ which is zero since $q$ is even. Thus 
$B = b_2(aE_{i_1j_1}+ E_{i_2j_2} + (1-a)E_{i_1j_2})$ has zero syndrome and hence $d(J_2,a)=3$ and part (b) holds.
Finally we may assume that, for each pair $(k,\ell)=(1,2), (2,3)$ or $(3,1)$, either $i_k=i_\ell$ or $j_k=j_\ell$.
Without loss of generality, and transposing if necessary, we may assume that $i_1=i_2$ (so $j_1\ne j_2$, since $B$ has weight three).
Suppose that $i_3\ne i_1$. Then we must have both $j_3=j_1$ and $j_3=j_2$, which is impossible. Hence also $i_3=i_1$, 
and $j_1,j_2,j_3$ are pairwise distinct.  However this implies that, for $i\ne i_1$, the entry $S_{ij_1}=b_1\ne0$. This contradiction completes the proof.  
\qed

{\prop Suppose $J_n\in\F_q^{n\times n}$ where $n\ge3$ or $q$ is odd,
and $a\in\F_q\setminus\{0,1\}$.
The parameters of $C(J_n,a)$ are $[n^2, (n-1)^2+1, 4]$ if
$\textup{char}(\F_q)\mid n$, and $[n^2,(n-1)^2,4]$ otherwise.}

\pf
Since $n\ge3$ or $q$ is odd, the minimum distance is $4$
by Theorem~\ref{thm:rank1}. However, $J_n^2=nJ_n$ implies that
$J_n$ is nilpotent if and only if the characteristic $p$ of $\F_q$ divides $n$.
Therefore by~\cite[Remark 2.10]{GP} the dimension of $C(J_n,a)$ is
$(n-1)^2+1$ if $p\mid n$, and it is $(n-1)^2$ otherwise.
\qed


\section{Automorphism group}
If $A$ and $A'$ are conjugate under the general linear group $\GL(n,F)$,
then the codes $C(A,a)$ and $C(A',a)$ have the same dimension, but almost
certainly different minimal distances. However, if $A$ and $A'$ are conjugate
by a monomial matrix, then $C(A,a)$ and $C(A',a)$ do have
the same minimal distances. The centralizer of $A$ in $\GL(n,F)$ induces
automorphisms of $C(A,a)$.

The \emph{adjacency matrix} of a graph $\Gamma$ with vertex set
$\{1,\dots,n\}$, is an $n\times n$ matrix with $(i,j)$-entry 1 or 0 according
as vertex $i$ and vertex $j$ are, or are not, adjacent in~$\Gamma$. 
Adjacency matrices can be interpreted as matrices over any field~$F$.

{\thm \label{QC} If $A\in F^{n\times n}$ is the adjacency matrix of a graph $\Gamma$ and $G=\textup{Aut}(\Gamma)$, then the direct product $G\times G$ acts
on the code $C(A,a)$ by coordinate permutations. }

\pf As is well-known \cite{B}, a permutation matrix $P$ lies in $\textup{Aut}(\Gamma)$ if and only if 
$PA=AP,$ that is if and only if $P\in C(A).$ Given $(P,Q)\in G\times G$
and $B\in C(A,a)$, it follows from Theorem \ref{nite}~(b) that
$P^{-1}BQ$ lies in $C(A,a).$ It is easy to verify that $B^{(P,Q)}=P^{-1}BQ$ 
defines an action of $G\times G$ on $C(A,a)$. This corresponds to the
so called `product action' of $G\times G$ permuting the $n^2$ coordinates
of the codewords of $C(A,a)$ (read off column by column). \qed

\medskip
A \emph{semiregular} permutation is a non-identity permutation, all cycles of which have the same length.
For a positive integer $n$ and a divisor $\ell$ of $n$, a code of length $n$ is called \emph{$\ell$-quasicyclic} if
a cyclic shift of each codeword by $\ell$ positions results in another codeword.  

{\cor If $A\in F^{n\times n}$ is the adjacency matrix of a graph admitting a semiregular automorphism with $m$ cycles, 
then $C(A,a)$ is up to equivalence $(n^2/m)$-quasicyclic.}

\pf
This follows from Theorem~\ref{QC} on taking the semiregular automorphism to act either on the rows or 
the columns of the matrix codewords.
\qed

{\cor If $A\in F^{n\times n}$ is a permutation matrix corresponding to a cycle of length $n$, then $C(A,a)$ is equivalent to an $n$-quasicyclic code.}

\pf Take $m=n$ in Corollary 4.2.
\qed

\medskip
Consider the following elementary observation. We call  the map $B \mapsto B^t$ the \emph{transposition permutation}. 
This permutation on $n^2$ elements can act on the coordinate entries of
$C(A,a)$ as the next result shows, \emph{c.f.} Theorem~\ref{nite}(e).

{\prop If $B \in C(A,a) $ and $a\neq 0,$ then $B^t\in C(A^t, a^{-1}).$ In particular if $A^t=\pm A$, then $C(A,1)$ and $C(A,-1)$ 
are invariant under the transposition permutation.}

\pf
Observe that $AB=aBA\Leftrightarrow A^tB^t=a^{-1} B^tA^t$, and $a=a^{-1}\Leftrightarrow a=\pm1$. 
\qed

\section{Examples}

Twisted centralizer codes give examples of many types of interesting
codes, including \emph{optimal} codes, see \cite[\S2.1]{HP}, and codes with
large minimal distance. Given $A\in F^{n\times n}$, the centralizer code $C(A,1)$
has minimal distance at most $n$ because $I_n\in C(A,1)$. By contrast, when
$a\ne0,1,$ the twisted centralizer codes $C(A,a)$ can have larger distances
as illustrated strikingly in Theorem~\ref{thm:dd}. This section provides
examples of codes $C(A,a)$ that are optimal and have $a\ne0,1$.
The results below were computed using \textsc{Magma}~\cite{BCP}
code~\cite{Glasby} and using Proposition~\ref{prop1.1}, where optimality
was confirmed using~\cite{G}.


\subsection{\texorpdfstring{$(-1)$}{}-centralizer codes over \texorpdfstring{$\F_3$}{}}
The codewords of $C(A,a)$ are easily described when $a=0,1$, see for
example~\cite{Stong} for $a=1$. In this subsection we consider twisted
centralizer codes over $\F_3$ or $\F_5$ with $a=-1$.

\bigskip
\noindent
{\bf Case 1: $n=3$.}
For each $A\in\F_3^{3\times 3}$ we compute the parameters $[N,k,d]$ for the
twisted centralizer codes $C(A,-1)$. Since $N=n^2=9$ we list the possible
pairs $[k,d]$ with the multiplicity as a superscript. The sum of the
multiplicities is $|\F_3^{3\times 3}|=3^9$.
\begin{align*}
  &[0,*]^{7722},[1,1]^{90},[1,2]^{720},[1,3]^{720},[1,4]^{900},
    [1,6]^{720},[1,9]^{360},[2,1]^{624},[2,2]^{1140},[2,3]^{480},\\
  &[2,4]^{1272},[2,5]^{384},[2,6]^{1248},[3,1]^{414},[3,2]^{876},
    [3,3]^{416},[3,4]^{840},[3,5]^{144},[3,6]^{40},\\
  & [4,1]^{216},[4,2]^{204},[4,4]^{48},[5,1]^{48},[5,2]^{48},[5,4]^{8},[9,1]^{1}. 
\end{align*}
The entry $[9,1]^{1}$ means that only $A=0$ has $C(A,-1)=\F_3^{3\times 3}$,
and $[0,*]^{7722}$ means that 7722 matrices $A\in\F_3^{3\times 3}$
have $C(A,-1)=\{0\}$. By convention, the minimal distance of the zero subspace of $F^{n\times n}$ is $d=n^2$. The above data were computed using
the computer code~\cite{Glasby}.

{\remk \rm The ternary codes with parameters $[9,5,4]$, $[9,3,6]$, $[9,2,6]$, $[9,1,9]$ are all optimal ternary codes according to \cite{G}. 
This contrasts with ordinary centralizer codes, where the minimum distances are at most $n$ (and here $n=3$).}

\bigskip
\noindent
{\bf Case 2: $n=4$.}
Suppose $A\in\F_3^{4\times 4}$. The codes $C(A_i,-1)$ where $A_i$ is
listed below are optimal ternary codes
$$
A_1 = \begin{pmatrix} 0&0&0&1\\0&0&2&0\\1&0&1&0\\1&2&0&1\end{pmatrix},\; A_2= \begin{pmatrix} 0 & 0 & 1 & 1 \\ 0 & 1
& 1 & 2 \\ 2 & 2 & 1 & 1 \\ 1& 0& 0 &2
\end{pmatrix},\; A_3=  \begin{pmatrix}0 & 1 & 1 & 1 \\ 2 & 0
& 1 & 2 \\ 2 & 1 & 2 & 0 \\ 1& 1& 0 &2
\end{pmatrix},\; A_4=  \begin{pmatrix}2 & 1 & 2 & 2 \\ 2 & 1
& 2 & 2 \\ 2 & 1 & 2 & 2 \\ 1& 2& 1 &1
\end{pmatrix},
$$
with parameters $[16,2,12]$, $[16,3,10]$, $[16,4,9]$ and $[16,10,4]$, respectively.

{\remk \rm Given a $4\times 4$ matrix $A$, the centralizer code $C(A)$ (with
$a=1$) can correct at most~1 error since $d(C(A))\le4$. By contrast,
$C(A_1,-1)$ above can correct 5 errors and the two codes have the same
information rates when $\dim C(A)=2$.}

\subsection{\texorpdfstring{$2$}{}-centralizer codes over \texorpdfstring{$\F_5$}{}} We give examples of optimal $\F_5$-linear codes of the form $C(A,2)$.

\bigskip
\noindent
{\bf Case 1: $n=2$.} If $A = \begin{pmatrix}1 & 1 \\ 4 & 4\end{pmatrix}
\in\F_5^{2\times 2}$, then
$C(A,2)$ is an $\F_5$-linear code with parameters $[4,2,3]$. Note that this code is also an MDS-code. Recall that it is 
impossible to have a one-error correcting code from the centralizer of a $2\times 2$ matrix, 
but we are able to do so with the $2$-twist.

\bigskip
\noindent
{\bf Case 2: $n=3$.} The matrices $A_1,A_2,A_3\in\F_5^{3\times 3}$ below
$$A_1= \begin{pmatrix} 0 & 1 & 1 \\ 1 & 1& 0 \\ 2 & 0 & 3\end{pmatrix},\quad
  A_2= \begin{pmatrix} 0 & 1 & 1 \\ 1 & 0& 2 \\ 2 & 1 & 0\end{pmatrix},\quad
  A_3= \begin{pmatrix} 1 & 1 & 1 \\ 1 & 1& 1 \\ 3 & 3 & 3\end{pmatrix},$$
give optimal codes $C(A_1,2), C(A_2,2), C(A_3,2)$ with parameters $[9,2,7], [9,3,6]$ and $[9,5,4]$, respectively.

\bigskip
\noindent
{\bf Case 3: $n=4$.} The matrices $A_1,A_2\in\F_5^{4\times 4}$ below
$$A_1= \begin{pmatrix}0&0&0&1\\ 0&0&4&0\\3&2&2&2\\4&3&4&4\end{pmatrix},\quad
  A_2= \begin{pmatrix}0&0&0&1\\0&0&4&0\\3&2&2&4\\3&3&1&2\end{pmatrix},
$$
give optimal codes $C(A_1,2), C(A_2,2)$ with parameters $[16,2,13], [16,3,12]$, respectively.

\section{Conclusion}

In this paper, we have introduced and studied twisted centralizer codes, a non-trivial generalization of the centralizer codes of \cite{centralizer}.
The incorrect dimension formula of \cite[Theorem~2.4]{centralizer} was replaced by lower and upper bounds on the dimension of $C(A,a)$ 
in \cite[Theorem~4.7]{GP} (Theorem 2.3).
The lower bound  is especially relevant when the spectrum of $A$ contains eigenvalues in the ratio of $a.$
These bounds were exploited in Section 2 to obtain explicit results for examples from Hadamard matrices. 
It would be worthwhile to
find more examples based on combinatorial matrices (adjacency matrices of graphs and designs). The absolute upper bound on the dimension of $C(A,1)$ 
\cite[Theorem 2.1]{centralizer} as a function of $n$ is not easy to generalize to $a \neq 1.$
The proof of the Kronecker
product expression for the parity-check matrix has been simplified. An upper bound on the minimum distance based on the concept of product codes has been derived.

Our computational evidence indicates that twisted centralizer codes can have
much higher minimal distances than centralizer codes. 
Thus they retain the computational advantages of centralizer codes while having much higher error-correction capacity. 
The error-correction itself can sometimes be more simply expressed as was demonstrated by the example of $J_n$ in Section~\ref{sub:rank1}.

{\bf Acknowledgement:} We would like to thank the referee for making helpful
suggestions. The fourth author thanks Prof. Rioul for helpful discussions. 

\end{document}